\documentclass[final,1p,times]{elsarticle}

\usepackage{amssymb}
\usepackage{multicol}
\usepackage{multirow}
\usepackage{array}
\usepackage{hyperref}
\usepackage{float}
\usepackage{caption}
\usepackage{tabularx}
\usepackage{subcaption}
\usepackage{booktabs}
\usepackage{color}

\usepackage{booktabs}
\usepackage{tabu}
\usepackage{threeparttable}
\usepackage{pifont}
\usepackage{graphicx}  
\usepackage{subfloat}

\usepackage{wrapfig}
\newcommand{\comment}[1]{}

\journal{Biomedical Signal Processing and Control}

\begin{document}

\begin{frontmatter}
\title{Comparison of One- Two- and Three- Dimensional CNN models for Drawing-Test-Based Diagnostics of the Parkinson's Disease}

\author[First]{Xuechao Wang \corref{cor1}}
\author[First]{Junqing Huang}
\author[First]{Marianna Chatzakou}
\author[Second]{Sven N{\~o}mm} 
\author[Second]{Elli Valla}
\author[Third]{Kadri Medijainen}
\author[Fourth,Fifth]{Pille Taba}
\author[Sixth]{Aaro Toomela}
\author[First,Seventh]{Michael Ruzhansky }

\address[First]{Department of Mathematics: Analysis, Logic and Discrete Mathematics, Ghent University, Ghent, Belgium (e-mail: {\{xuechao.wang, junqing.huang, marianna.chatzakou, michael.ruzhansky\}@ ugent.be).}}
\address[Second]{Department of Software Science, School of Information Technologies, Tallinn University of Technology, Akadeemia tee 15a, 12618, Tallinn, Estonia (e-mail: \{elli.valla, sven.nomm\}@ taltech.ee).}
\address[Third]{Institute of Sport Sciences and Physiotherapy, University of Tartu, Puusepa 8, Tartu 51014, Estonia (e-mail: kadri.medijainen@ut.ee)}
\address[Fourth]{Department of Neurology and Neurosurgery, University of Tartu, Puusepa 8, Tartu 51014, Estonia (e-mail: pille.taba@kliinikum.ee)}
\address[Fifth]{Neurology Clinic, Tartu University Hospital, Puusepa 8, Tartu 51014, Estonia}
\address[Sixth]{School of Natural Sciences and Health Tallinn University, Narva mnt. 25, 10120 (e-mail: aaro.toomela@tlu.ee)}
\address[Seventh]{School of Mathematical Sciences, Queen Mary University of London, London, UK (e-mail: m.ruzhansky@qmul.ac.uk)}

\cortext[cor1]{Corresponding author}



\begin{abstract}
\noindent{\em Subject:} In this article, convolutional networks of one, two, and three dimensions are compared with respect to their ability to distinguish between the drawing tests produced by Parkinson's disease patients and healthy control subjects.

\noindent{\em Motivation:} The application of deep learning techniques for the analysis of drawing tests to support the diagnosis of Parkinson's disease has become a growing trend in the area of Artificial Intelligence. 

\noindent{\em Method:} The dynamic features of the handwriting signal are embedded in the static test data to generate one-dimensional time series, two-dimensional RGB images and three-dimensional voxelized point clouds, and then one-, two-, and three-dimensional CNN can be used to automatically extract features for effective diagnosis. 

\noindent{\em Novelty:} While there are many results that describe the application of two-dimensional convolutional models to the problem, to the best knowledge of the authors, there are no results based on the application of three-dimensional models and very few using one-dimensional models.

\noindent{\em Main result:} The accuracy of the one-, two- and three-dimensional CNN models was $62.50\%$, $77.78\%$ and $83.34\%$ in the {\em DraWritePD} dataset (acquired by the authors) and $73.33\%$, $80.00\%$ and $86.67\%$ in the {\em PaHaW} dataset (well known from the literature), respectively. For these two data sets, the proposed three-dimensional convolutional classification method exhibits the best diagnostic performance.
\end{abstract}





\end{frontmatter}


\section{Introduction}\label{sec:introduction}
The present paper compares the one-, two- and three-dimensional deep convolutional neural network (CNN) models for the analysis of the drawing tests used to support the diagnosis of Parkinson's disease (PD).  Parkinson's disease is one of the most common neurodegenerative disorders. Its symptoms like rigidity, tremor, and non-purposeful motions severely affect the quality of patient's life. Although at the time of writing of this article there is no cure for PD, proper therapy may allow one to eliminate these symptoms or reduce their effect on motion and improve the quality of daily life. Drawing tests have been used to diagnose PD and assess its severity for nearly a century \cite{phillips1991can}. These tests require one to continue, copy, or trace the repeating pattern or contour of an object. Only paper and pen were required to perform the test, while the evaluation was performed visually by the practitioner. This method is limited by the practitioner's experience, the ability of the naked eye, and the fact that the smoothness parameters of the drawing could not be recorded for future analysis and comparison. The seminal paper \cite{marquardt} laid the basis for computer-aided analysis of drawing and writing tests. In \cite{marquardt} it was suggested to use a digital table to acquire time-stamped coordinates of the tip of the stylus and compute the kinematic parameters that describe the movements of the tip of the stylus relative to the device screen. After that, kinematic parameters can be computed that describe the drawing movements observed during the test. In addition to providing kinematic and pressure descriptions of the testing process, digitisation of the tests offers the possibility of performing testing before the visit to the doctor, saving valuable time and providing access for medical professionals. The digitisation of the testing procedures has greatly expanded the set of features \cite{drotar} and demonstrated the importance of features based on tremor\cite{VALLA2022103551} to achieve highly accurate results. The results mentioned above use statistical machine learning techniques, whereas the data are presented in tabular form. Without undermining the importance of these results, it is important to mention that the evaluation procedure is different from that used by the human practitioner. This makes it difficult for the practitioner to interpret the results of digitised tests. The core difference is that a human mostly assesses the shapes of the drawn contours, smoothness of the movements (as the naked eye can see), and presence or absence of the errors, whereas statistical learning algorithms use the set of values describing the kinematic and pressure parameters. One way to mimic the human practitioner is to employ deep learning techniques to classify drawings according to their shapes \cite{moetesum2019assessing}. Such a solution would be closer to human assessment. However, it would ignore the advantages of performing the test on digital tables or tablet PCs that can capture the kinematic and pressure parameters that describe the test. The colouration of the drawn lines was proposed in \cite{rios2019analysis} to encode the pressure parameter. This may be seen as the bridge between mimicking analysis made by human practitioners and novel approaches that are based on features described by kinematic and pressure properties of the motion.  Later, \cite{NOMM2020260}, \cite{zarembo2021cnn} expanded on this idea and suggested varying the thickness of the drawn contours to encode one more kinematic or pressure parameter in the drawings. Later in \cite{galaz2022comparison}, the ``hand-crafted'' and CNN-learnt features were compared. These approaches assume that the data are provided in the form of images and that the CNN classifier is used to estimate whether the test drawing was produced by the PD patient or the healthy control (HC) subject. In theory, encoding more kinematic and pressure parameters in the drawing could further increase the goodness of the diagnostic support model. Following this idea, applications of three-dimensional convolutional neural networks seem to be a logical step. At the same time, one-dimensional CNN has been successfully applied to similar problems \cite{gil2019parkinson}, leading to the idea of comparing the three cases. The CNN structures best suited for the particularities of the images resulting from the drawing testing procedures will be selected first. Then, different feature sets will be selected to encode in the original drawings. The selected models will then be trained and validated to determine the quality of the models and the required training time. The organisation of the paper follows the classical academic style. An overview of the literature necessary to position the current contribution and explain its novelty is presented in Section \ref{sec:literature}. Background information explaining how symptoms of PD influence the feature engineering process is presented, together with a description of the hardware and software settings used for data acquisition in Section \ref{sec:experimental}. The choice of CNN architectures to compare with all other parameters of computational experiments is presented in Section \ref{sec:workflow}. The main results are stated in Section \ref{sec:results}. The results achieved and their medical interpretation are discussed in Section \ref{sec:discussion}. Conclusions are drawn in the last section.    

\begin{table*}[htbp]
	\centering
	\footnotesize
    \caption{Overview of the related works.} 
    \label{tab:related works}
		\begin{tabular}{l l l l l l}
			\toprule
			Author(s) & Dataset & Features & Models & Accuracy & Year \\
            \hline
               Drotár et al. \cite{drotar2014analysis} & PaHaW & kinematic features & SVM & $85.61$ & $2014$ \\
               Drotár et al. \cite{drotar2014decision} & PaHaW & kinematic and & SVM  & $88.13$ & $2015$\\
               &  & spatio-temporal features &  &  &\\
              Drotár et al. \cite{drotar} & PaHaW & kinematic and pressure & SVM & $87.4$ & $2016$ \\
               Pereira et al. \cite{pereira2016new} & HandPD & kinematic features & SVM & $65.88$ & $2016$ \\
               Pereira et al. \cite{pereira2018handwritten} & HandPD & time series based features & 2D CNN & $96.35$ & $2018$ \\
               Gil-Martín et al. \cite{gil2019parkinson} & \cite{isenkul2014improved} & kinematic features & 1D CNN & $96.5$ & $2019$ \\ 
              Diaz et al. \cite{diaz2019dynamically} & PaHaW & dynamically enhaced static & 2D CNN + SVM & $75$ & $2019$ \\
              Naseer et al. \cite{naseer2020refining} & PaHaW & fine-tuned-ImageNet features & AlexNet & $98.28$ & $2020$ \\
              Diaz et al. \cite{diaz2021sequence} & PaHaW & raw and derived features & 1D CNN-BiGRU & $93.75$ & $2021$ \\
               Gazda et al. \cite{gazda2021multiple} & PaHaW & fine-tuned-ImageNet features & 2D CNN & $85.7$ & $2022$ \\
               & HandPD &  &  & $92.7$ & \\
			\toprule
		\end{tabular} 
\end{table*}

\section{Literature overview and state of the art}\label{sec:literature}

The CNN concept was formalised in \cite{IEEE98_CNN_Lecun} for two-dimensional image recognition. Now, nearly $30$ years later, there are numerous types of CNN architecture \cite{esteva2021deep} used for image recognition and video processing. Although most CNN types are within two-dimensional CNN models, one can distinguish one-dimensional and three-dimensional CNN types \cite{alzubaidi2021review}.

Currently, the most popular approach to investigate the potential of automated handwriting analysis for the diagnosis of PD involves the use of dynamic information in the handwriting process to produce a more discriminating feature set of different dimensional data.



Dynamic handwriting analysis benefits from the use of digital tablets and electronic pens \cite{pereira2019survey}. Using these devices, it is possible to directly measure the temporal and spatial variables of handwriting, the pressure applied to the writing surface, the inclination of the pen, and the movement of the pen when it is not in contact with the surface (i.e. in air) \cite{drotar2014analysis}. The applicability of kinematic, geometric, and non-linear dynamic characteristics was explored in a model of handwriting impairment in PD patients \cite{rios2019analysis}, and dynamic pressure features were encoded as colour. Furthermore, aerial movement during handwriting has a significant impact on the precision of disease classification \cite{drotar2014analysis}.

On the other hand, since dynamic analysis needs to consider not only the underlying generative process but also the geometry of handwritten patterns, graph analysis in a two-dimensional space is usually preferred to signal analysis in a one-dimensional space. 

  


Encouraged results have recently been reported to quantitatively assess the visual properties of handwritten motion samples from patients with PD using raw, filtered median and edge images\cite{moetesum2019assessing}. However, according to others \cite{galbally2015line} a better understanding can be obtained using ``dynamic augmentation'' of static handwriting. Instead of simply using images of handwritten patterns, less realistic but more discriminating images are obtained by including additional dynamic information in the generation process.

One of the most serious problems to solve before applying deep learning techniques in the analysis of drawing and writing tests is the small size of the data sets available for training and validation. Due to the differences between testing protocols used in different medical centres and strict data handling requirements, acquiring sufficiently large data sets for training is not a viable solution. On the side of deep learning, there are two techniques that are used to overcome this problem. The first technique is data augmentation \cite{shorten2019survey}. The augmentation procedure is based on the application of affine transformations, local nonlinear distortions, colour alternation, and noising to each image of the data set many times, whereas each alternated clone inherits the label of the original image. The set of transformations and their magnitude is chosen at random. This method was used in \cite{NOMM2020260} and \cite{zarembo2021cnn} and \cite{kamran2021handwriting}. The latest has provided an analysis of different transformation types and their influence on modelling quality. Alternatively to this unsupervised technique, applications of generative adversary networks (GANs) \cite{DZOTSENIDZE2022108,chen2022generative} may be used. Table \ref{tab:related works} summarises the main characteristics of previous works on PD diagnosis based on digital drawings: reference, dataset, feature set, method, performance, and year.

\begin{figure}[htbp]
	\centering  
	\includegraphics[width=0.9\textwidth, height =0.15\textheight]{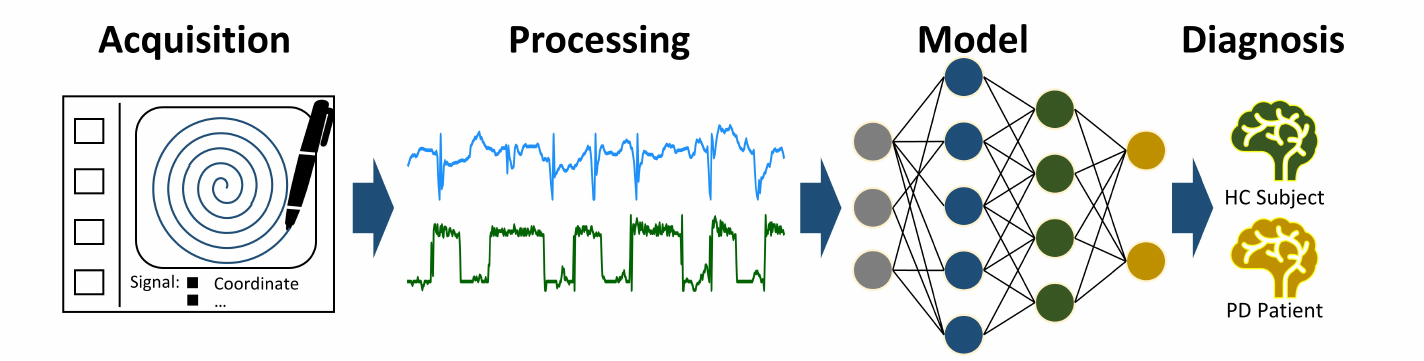}
	\caption{The workflow for diagnosis of Parkinson's disease.}  
	\label{fig:workflow}
\end{figure}

\section{Problem statement}
Main motivation of the present research is that, while 2D CNN is the dominant type when talking about image analysis, 3D models may allow for the encode of more kinematic and pressure parameters of the motion. Of course, one may suggest increasing the dimensionality even higher and using all the available parameters, but higher-dimensional convolutions are difficult to interpret. Therefore, CNNs with dimensions larger than three are left outside of the current research framework. On the contrary, one-dimensional CNNs have been used successfully before and are easy to interpret. This leads to the formal problem statement of the present investigation. Compare the performance of one-, two-, and three-dimensional CNN models for spiral drawing tests classification. This requires one to answer the following research questions. 
\begin{enumerate}
    \item Choose data set enhancement technique to encode different kinematic and drawing parameters into the drawing.
    \item Choose the data set augmentation technique such that it acts in a similar way for one-, two-, and three-dimensional cases.
    \item Choose the feature set(s) to encode. 
    \item Choose the CNN models structures which are the most similar among the one-, two-, and three-dimensional cases.
\end{enumerate}

\section{Materials}\label{sec:experimental}
In this work, two data sets have been considered. The first data set, called {\em DraWritePD} \cite{valla2022tremor}, was acquired by the authors. The second data set, known as {\em PaHaW}, was kindly provided by the authors of \cite{drotar2014analysis, drotar2016evaluation}. Both datasets use similar digital signal acquisition equipment, and the handwriting data contains the same dynamic features (time sequences). In addition, they all contain a similar number of samples from each class, making the experiments more balanced.

\subsection{DraWritePD} 
The ``Drawing and handwriting tests for Parkinson’s diagnostics'' ({\em DraWritePD}) collects handwriting data from $25$ patients with PD and $34$ healthy control (HC) subjects of the same age and sex. For the group of patients with PD, the mean age was $74.1$ $\pm$ $6.7$ years. For the group of subjects with HC, the mean age was $74.1$ $\pm$ $9.1$ years. To acquire handwriting signals, special applications were designed for the Apple IPad Pro ($9.6$ inch, $2016$ year) with the first generation of the Apple Pen. The application displays test instructions and the reference drawing on the IPad screen and records dynamic information from the Apple Pen tip accompanied by the time stamp. PD patients and their HC counterparts were asked to complete a series of handwriting tests consisting of $12$ different drawing and writing tasks. When the task was completed, this dynamic information was stored in the file for later processing. It may be seen as a matrix with rows corresponding to the timestamps and columns corresponding to independent dynamic features, including: x coordinate (mm); y coordinate (mm); timestamp (sec); pressure (arbitrary unit of force applied on the surface: [0,$\cdots$, 6.0]); altitude (rad); azimuth (rad). In the present research, only digital versions of the Archimedes spiral drawing test (ASD) were considered.  

\subsection{PaHaW} 
The ``Parkinson’s disease handwriting database'' ({\em PaHaW}) collects handwriting data from $37$ patients with PD and $38$ HC subjects \cite{drotar2014analysis, drotar2016evaluation}. No significant differences were found between the groups with respect to age or sex. The database was acquired in cooperation with the Movement Disorders Centre of the First Department of Neurology, Masaryk University, and St. Anne's University Hospital in Brno, Czech Republic. Each subject was asked to complete multiple handwriting tasks according to the prepared filled template at a comfortable speed. A tablet was overlaid with an empty paper template (containing only printed lines and a square box specifying the area for the Archimedean spiral), and a conventional ink pen was held in a normal fashion, allowing for immediate full visual feedback. Handwriting signals were recorded using an Intuos $4$M (Wacom technology) digitising tablet at a sampling frequency of $150$ Hz during pressure on the writing surface and movement over the writing surface. We denote these signals by on-surface movement and on-air movement, respectively. The recordings started when the pen touched the surface of the digitiser and finished when the task was completed. The tablet captured the following independent dynamic features: x coordinate; y coordinate; timestamp; button status; pressure; altitude; and azimuth. The button status was a binary variable, being $0$ for pen-up state (in-air movement) and $1$ for the pen-down state (on-surface movement). Although the task set presented in the {\em PaHaW} dataset is quite different from that used in the {\em DraWritePD} dataset, ASD was present in both datasets and was therefore used in this work.

\section{Research workflow}\label{sec:workflow}
In this section, we continue the idea of embedding dynamic handwriting features into static handwriting images. Specifically, more kinematic and pressure features are gradually encoded to generate higher-dimensional data representations, so that the corresponding one-, two-, and three-dimensional convolutional neural networks are used to analyse handwriting tasks to support the diagnosis of PD. Figure \ref{fig:workflow} shows a general overview of the proposed automatic PD diagnosis system. Details of each stage are presented in the following subsections.

\begin{figure}[htbp]
    \centering
    \subfloat[$1$D - Time Series]{
        \includegraphics[width=0.31\textwidth]{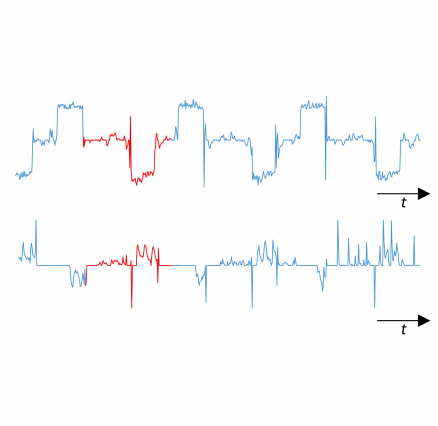}}
    \subfloat[$2$D - RGB Image]{
        \includegraphics[width=0.31\textwidth]{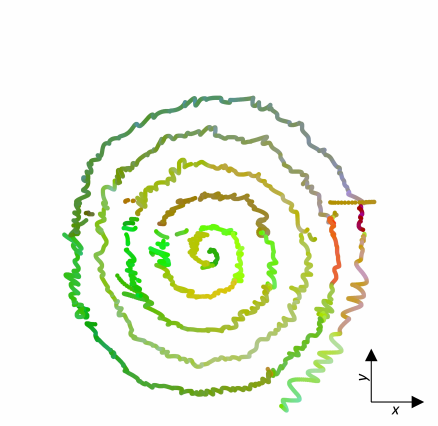}}
    \subfloat[$3$D - Point Cloud]{
        \includegraphics[width=0.31\textwidth]{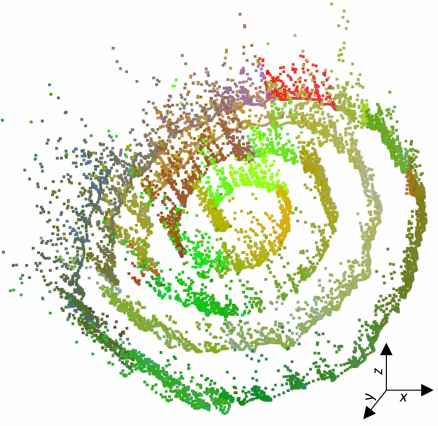}}
    \caption{Enhanced data in different dimensional cases. In the one-dimensional ($1$D) case, the handwriting signal is enhanced into a time series, in which the raw dynamic features (such as x-coordinate and y-coordinate) are directly used; in the two-dimensional ($2$D) case, the handwriting signal is enhanced into an RGB image, in which the coordinate features (x-coordinate, y-coordinate) are used as (x,y) position information, and the (azimuth,altitude,pressure) features are used as (R,G,B) colour information, and the velocity feature is used as line width information; in the three-dimensional ($3$D) case, the handwriting signal is enhanced into a point cloud, in which the features (x-coordinate, y-coordinate, velocity) are used as (x,y,z) position information, and the (azimuth,altitude,pressure) features are used as (R,G,B) colour information.}
	\label{fig:data enhancemnet}
\end{figure} 

\subsection{Data Processing}\label{sec:data processing}
Data processing consists of three main steps: data preparation, data enhancement, and data augmentation. Note that the following are introduced through the $1$D, $2$D, and $3$D cases, respectively. First, since the raw dataset inevitably contains some features that are not suitable for direct use, such as the units of the two datasets being different, the raw dynamic features are preprocessed with maximum and minimum normalisation before data enhancement to convert them into the same unit. Subsequently, the data enhancement gradually encodes dynamic features to generate enhanced data of different dimensions. Specifically, for the $1$D case, the data encoding method is to directly regard the raw dynamic features (such as the x coordinate and the y coordinate) in the handwriting signal as $1$D time series data. However, it should be noted that the encoding method of the timestamp feature is to replace the timestamp feature itself with the velocity feature calculated by combining the timestamp feature and the coordinate features. For the $2$D RGB image encoding method, the coordinate features (the x coordinate and y coordinate) used in the $1$D case are used as the pixel position information corresponding to each data point. Moreover, the azimuth, altitude, and pressure features of each data point are used as the red (R), green (G), and blue (B) colour information of the corresponding pixel. The velocity feature is encoded as line width information. For the $3$D case, the only difference from the $2$D case is the location information for each data point. It not only utilises the coordinate features (the x coordinate and y coordinate), but also combines the time feature (timestamp) to calculate the velocity feature of each data point, thus adopting (x coordinate, y coordinate, velocity) as the $3$D position information of each data point. Afterwards, the generated raw point cloud data is voxelized into a matrix form acceptable to the CNN model with a fixed grid resolution (for convenience, hereinafter referred to as point cloud). It is worth noting that CNN has a powerful feature extraction ability \cite{alzubaidi2021review}, so the raw dynamic features were used directly for the enhancement of the data and no additional hand-crafted features were designed except for the velocity feature. Figure \ref{fig:data enhancemnet} shows the results of the data enhancement in different dimensions. Furthermore, a major challenge for the diagnosis of PD is the lack of suitable data. The direct application of CNN cannot effectively process raw handwriting signals collected from patients. One current approach to address this challenge is to augment data through data augmentation techniques or by combining multiple datasets \cite{kamran2021handwriting}, or employ pre-trained transfer learning strategies \cite{naseer2020refining}. In the present work, data augmentation techniques are employed to significantly increase the diversity of PD handwriting samples, which can be roughly classified into the following categories; original, flipping, rotation, illumination, and jitter.
\begin{itemize}
\item Flipping: Flipping produces a mirror data, where the RGB image is flipped horizontally and vertically, and the point cloud is flipped along the x-axis and y-axis, respectively.
\item Rotation: The data are rotated by a given angle, such as 90, 180, or 270 degrees, where the RGB image is rotated around the centre point, and the point cloud is rotated around the z-axis.
\item Illumination: Illumination can be implemented by adjusting the colour map (RGB values), where random values are added to the R, G, and B channels of RGB images and point clouds.
\item Jitter: Jitter can be implemented by adjusting the values of the coordinate features in the handwriting signal, where random values are added to the values of the coordinate features, and the resulting signal is enhanced.
\end{itemize}
Note that the first three data augmentation techniques are not suitable for the $1$D case. Furthermore, relatively large data can negatively affect training time, while concise features can lead to under-fitting. Based on the available data, we first propose that the data sizes of $1$D, $2$D, and $3$D are resized to $128$, $128^2$, and $128^3$, respectively.

\begin{figure}[htbp]
    \centering
    \subfloat[$1$D - Conv1d Kernel]{
        \includegraphics[width=0.31\textwidth]{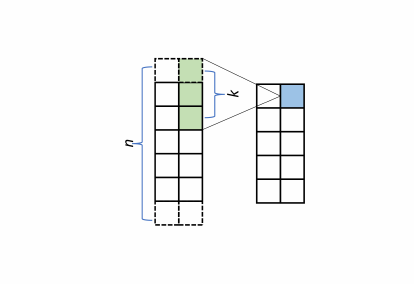}}
    \subfloat[$2$D - Conv2d Kernel]{
        \includegraphics[width=0.31\textwidth]{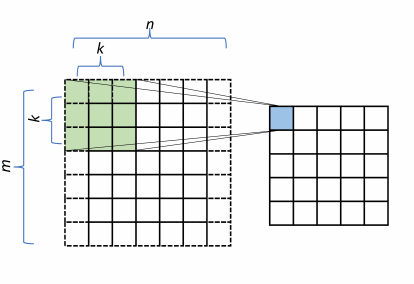}}
    \subfloat[$3$D - Conv3d Kernel]{
        \includegraphics[width=0.31\textwidth]{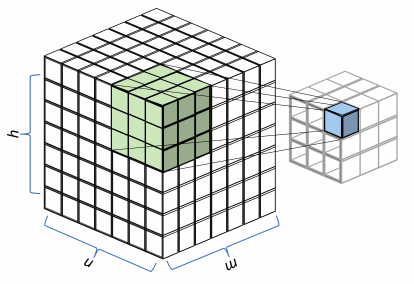}}
    \caption{The schematic diagram of convolution operation in different dimensions. The green area represents the convolution area, and the blue area represents the convolution result, where $m$, $n$, and $h$ are the feature map sizes, and $k$ represents the convolution kernel size.}
	\label{fig:kernel size}
\end{figure} 

\subsection{Neural Network}
The convolutional neural networks are bioinspired variants of multilayer perceptrons (MLP) that can perform a variety of machine learning tasks without requiring the user to design and provide any hand-crafted features \cite{alzubaidi2021review}. Recently, due to the development of new CNN variants, they have shown promising performance in traditionally challenging tasks with breakthrough progress \cite{wu2016comprehensive,miao2018dilated}. The paradigm-shifting results provided by CNN are done in part with the help of extremely large training datasets. However, as mentioned earlier, one of the biggest limitations in the medical community is the inability to access larger, labelled, high-quality data that are sensitive, confidential, and difficult to collect. Due to the insufficient amount of data, in this work, in addition to using data augmentation techniques to augment data, we also incorporate a simplified version of the AlexNet \cite{krizhevsky2017imagenet} architecture, which consists of two main parts (convolutional layers for feature extraction and fully connected layers for classification), similar to \cite{naseer2020refining}. Furthermore, it is worth pointing out that for fair comparisons, the one-, two-, and three-dimensional convolutional neural networks use the same network architecture, and the only difference is the size of the convolution kernel and the convolution method. 


\begin{figure}[htbp]
	\centering  
	\includegraphics[width=0.99\textwidth, height =0.2\textheight]{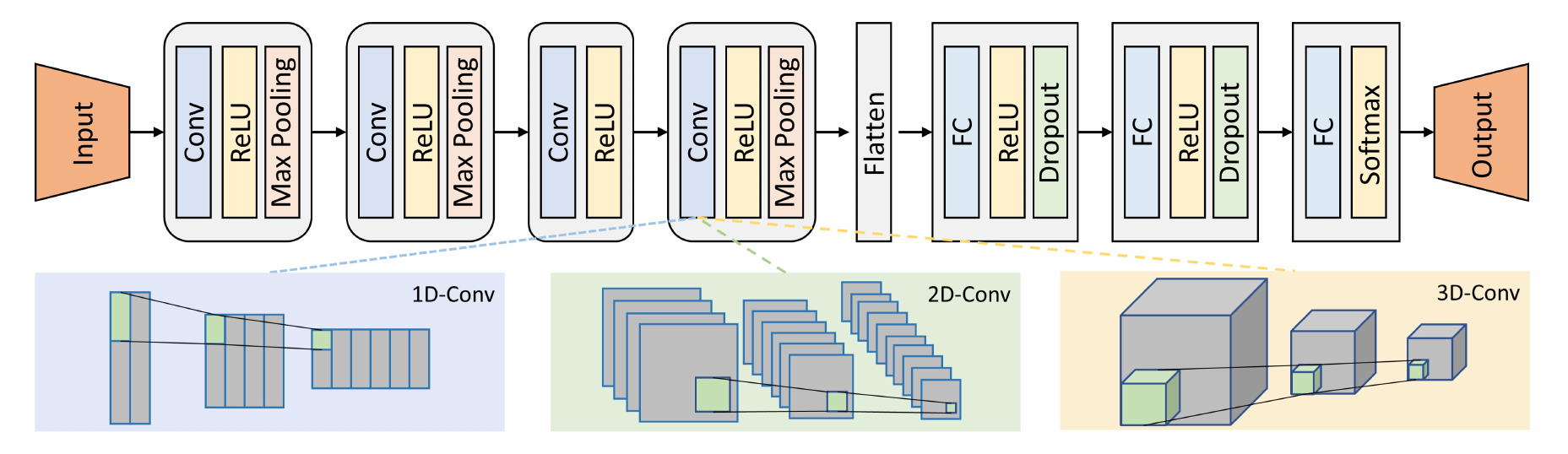}
	\caption{The convolutional neural network model, in which the one-, two-, and three-dimensional convolutional neural networks use the same model architecture, only the convolution method is different.}  
	\label{fig:model}
\end{figure}

\begin{table}[htbp]
	\footnotesize
    \centering
    \begin{threeparttable}
    \caption{Architectural differences between one-, two-, and three-dimensional convolutional neural networks.}
    \label{tab:model}
    \begin{tabular}{l c c l l l l l l l l} 
		\toprule  
        & & & \multicolumn{2}{l}{$1$D CNN} & & \multicolumn{2}{l}{$2$D CNN} & &\multicolumn{2}{l}{$3$D CNN} \\
        \cline{4-5} \cline{7-8} \cline{10-11}
        Layers & Filter & S & Input & K & & Input & K & & Input & K\\
        \hline
        Conv+ReLU  & $48$  & $2$ & $(6,128)$ & $5$ & & $(3,128,128)$ & $(5,5)$ &  & $(3,128,128,128)$ & $(5,5,5)$ \\
        MaxPooling  & - & - & $(48,64)$ & - && $(48,64,64)$ &-& & $(48,64,64,64)$ &- \\
        Conv+ReLU  & $128$  & $2$ & $(48,32)$ & $5$ & & $(48,32,32)$ & $(5,5)$ & & $(48,32,32,32)$ & $(5,5,5)$ \\
        MaxPooling  & - & - & $(128,16)$ & - & & $(128,16,16)$ & - & & $(128,16,16,16)$ &- \\
        Conv+ReLU  & $192$  & $1$ & $(128,8)$ & $3$ & & $(128,8,8)$ &$(3,3)$&& $(128,8,8,8)$&$(3,3,3)$ \\
        Conv+ReLU  & $192$  & $1$  & $(192,8)$ & $3$ & & $(192,8,8)$ &$(3,3)$&& $(192,8,8,8)$ & $(3,3,3)$ \\
        MaxPooling  &  - & - &$(192,8)$&-& & $(192,8,8)$ &-& & $(192,8,8,8)$&-\\
        Flatten & -& - & $(192,4)$&- && $(192,4,4)$&-& & $(192,4,4,4)$&- \\ 
        FC+ReLU &  - & - & $768$&- && $3072$&-& & $12288$&- \\
        Dropout &  -& - & $192$ &-&& $192$&-& & $192$&- \\
        FC+ReLU &  -& - & $192$ &-&& $192$&-& & $192$&-\\ 
        Dropout &  -& - & $128$ &-&& $128$&-& & $128$&- \\
        FC+Softmax &  -&- & $128$ &-&& $128$&-& & $128$&- \\
         \hline
         Params &  &  & $0.39$MB & & & $1.33$MB && & $4.83$MB&\\
        \toprule 
	\end{tabular}
    \begin{tablenotes}
        \footnotesize
        \item \textit{Note}: The abbreviations Conv, ReLU, and FC denote the convolutional layer, Rectified Linear Unit, and fully connected layer; and K and S are the kernel size and stride size; Params $=$ parameters .
	\end{tablenotes}
    \end{threeparttable}
\end{table}
The simplified AlexNet architecture consists of four convolutional layers, a maximum pooling layer, a dropout layer, and three fully connected layers. The output of the fully connected layer is passed to a softmax layer to produce a distribution on the labels of the $2$ class. The first convolutional layer uses $48$ kernels of size $5$ with a stride of $2$. The second convolutional layer takes as input the output of the first layer (ReLU activation and Max pooling) and filters it using $128$ kernels of size $5$. The third layer has $192$ kernels of size $3$ connected to the activated and pooled outputs of the second layer, while the fourth layer contains $192$ kernels of size $3$. The dropout layer in the fully connected layer temporarily removes nodes from the network with probability $50\%$ during network training. The specific convolution kernel and the convolution operation are shown in Figure \ref{fig:kernel size}. Details of the simplified AlexNet architecture deployed in our experiments are shown in Table \ref{tab:model} and Figure \ref{fig:model}. 

\subsection{Implementation}
Model training and testing were carried out on a PC that has an Intel(R) Core(TM) i$7$-$11700$K CPU with $3.60$ GHZ($8$ CPU), $32$GB RAM and an NVIDIA GEFORCE RTX$3070$Ti graphics card with $8$ GB memory. The initial learning rate was set to $1e-4$, and Adam \cite{kingma2014adam} optimiser was used to train the model. The cross-entropy loss function was used to optimise the model parameters.  $80\%$ of the original data set is for training and validation, and $20\%$ is for testing. Various metrics were used to measure the performance of the model in more detail, including accuracy, precision, sensitivity, specificity, and $F_{1}$-score (see \cite{grandini2020metrics}).

\section{Main results}\label{sec:results}

\begin{figure}[htbp]
    \centering
    \subfloat[$1$D - DraWritePD]{
        \includegraphics[width=0.45\textwidth]{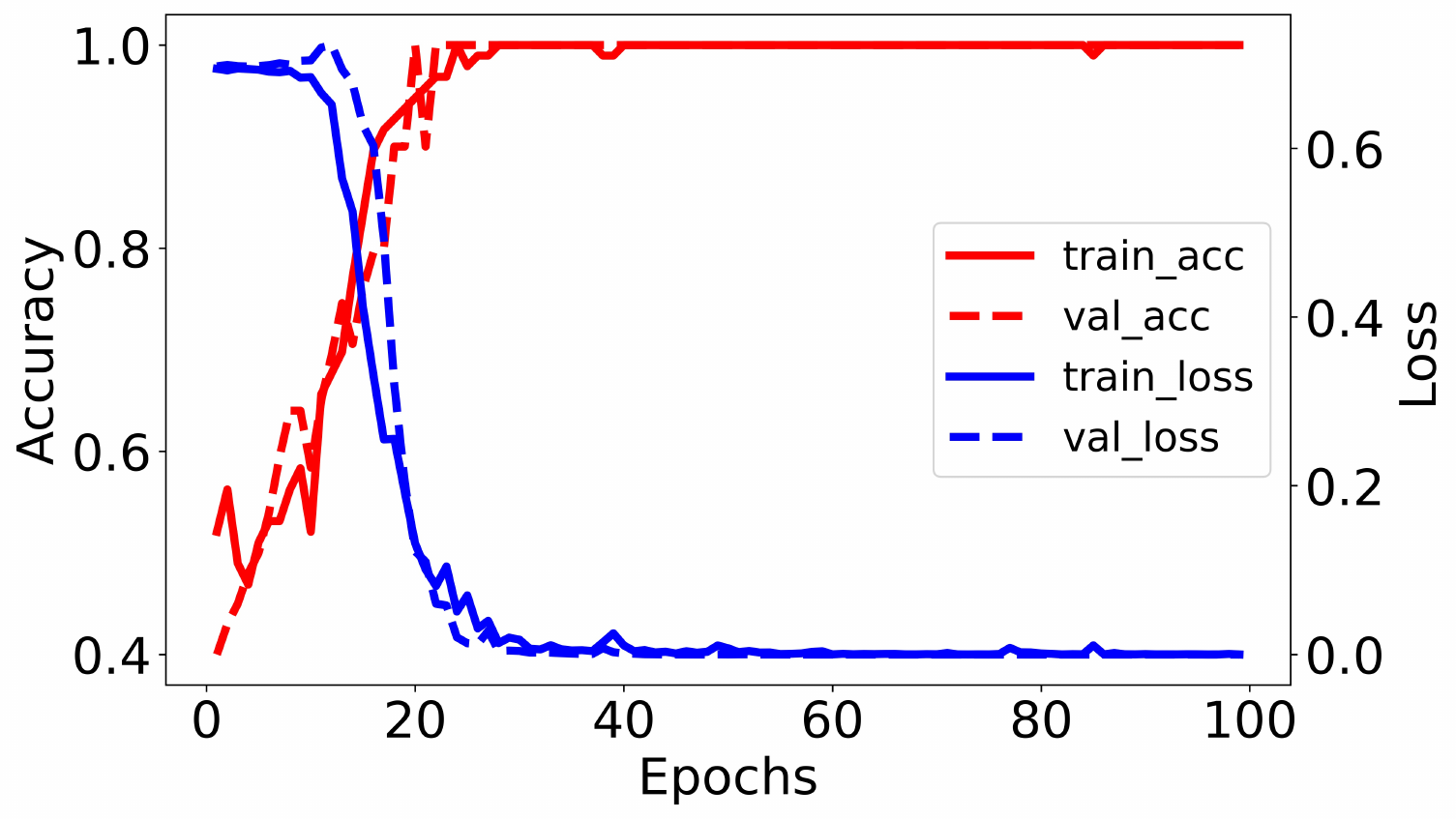}}
    \subfloat[$1$D - PaHaW]{
        \includegraphics[width=0.45\textwidth]{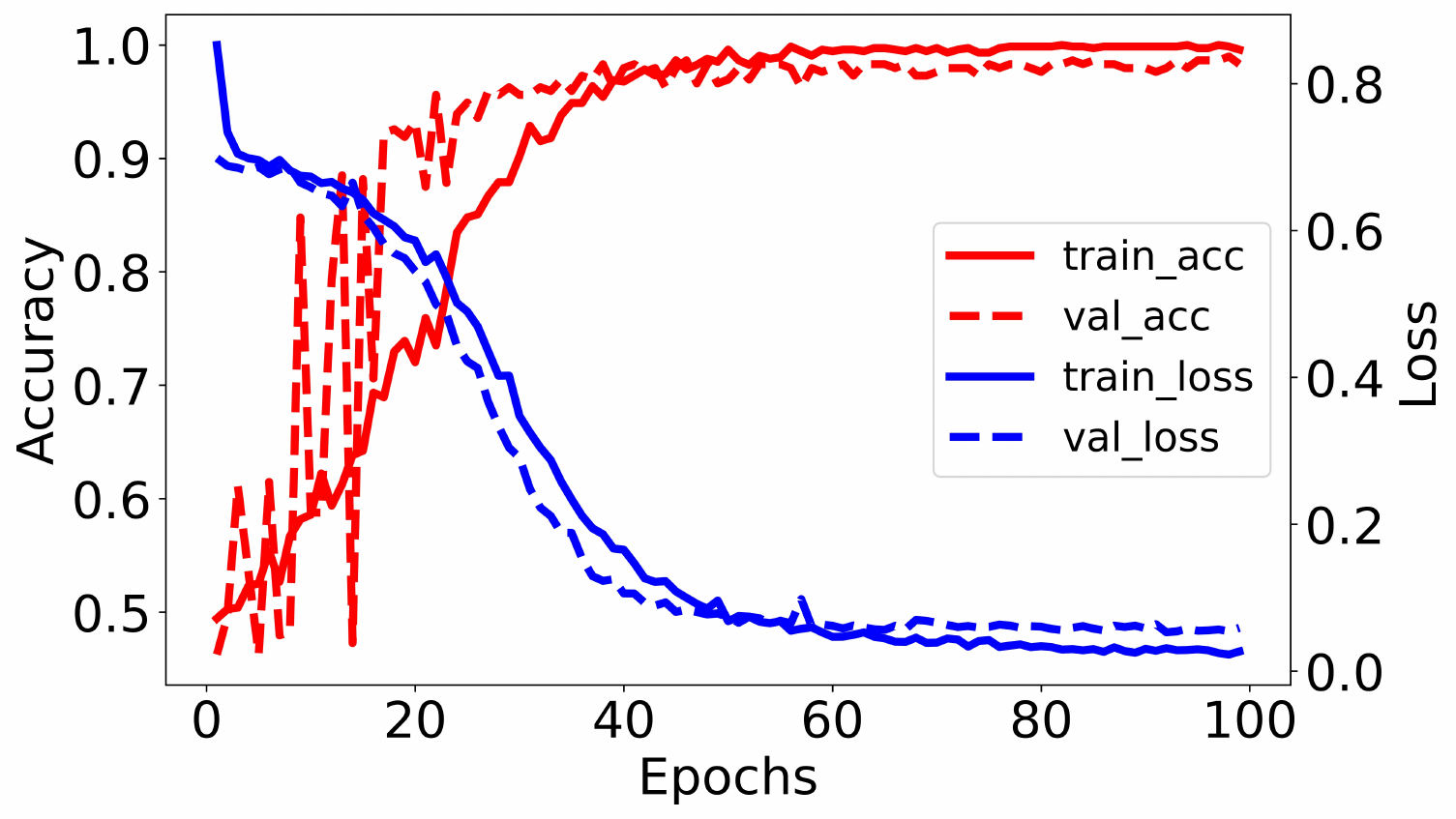}}
        \\
    \subfloat[$2$D - DraWritePD]{
        \includegraphics[width=0.45\textwidth]{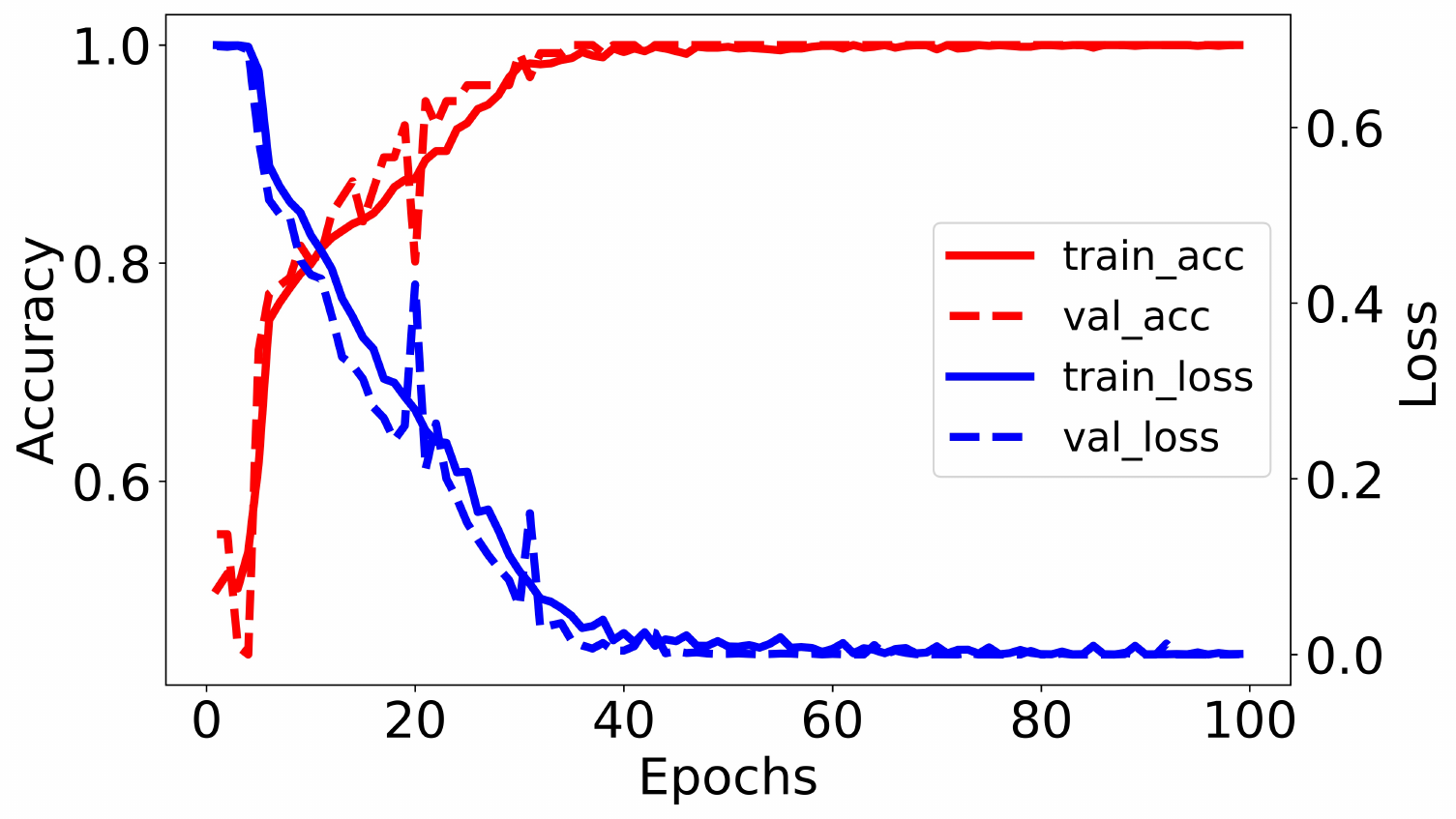}}
     \subfloat[$2$D - PaHaW]{
        \includegraphics[width=0.45\textwidth]{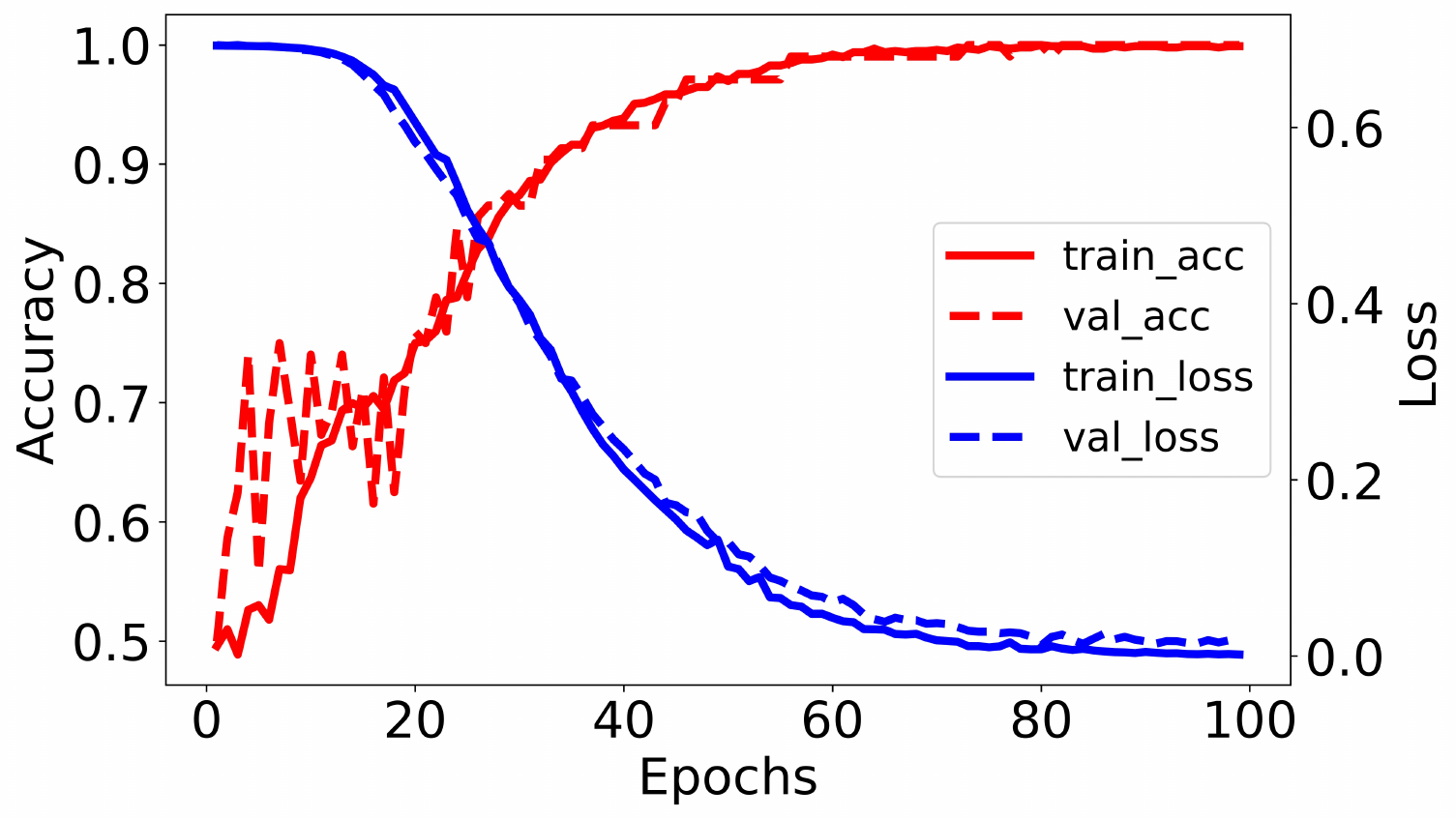}}
        \\
    \subfloat[$3$D - DraWritePD]{
        \includegraphics[width=0.45\textwidth]{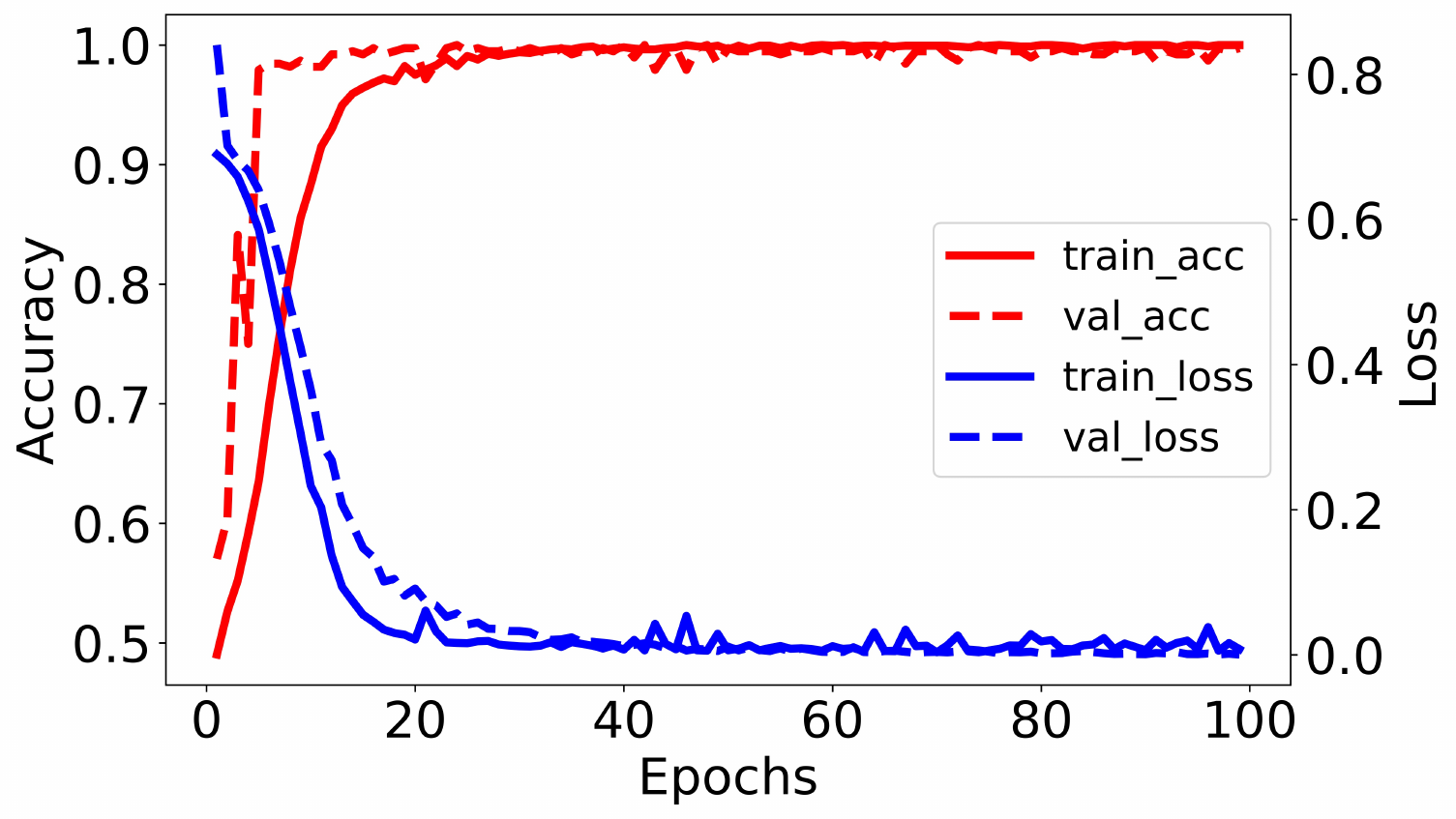}}    
    \subfloat[$3$D - PaHaW]{
        \includegraphics[width=0.45\textwidth]{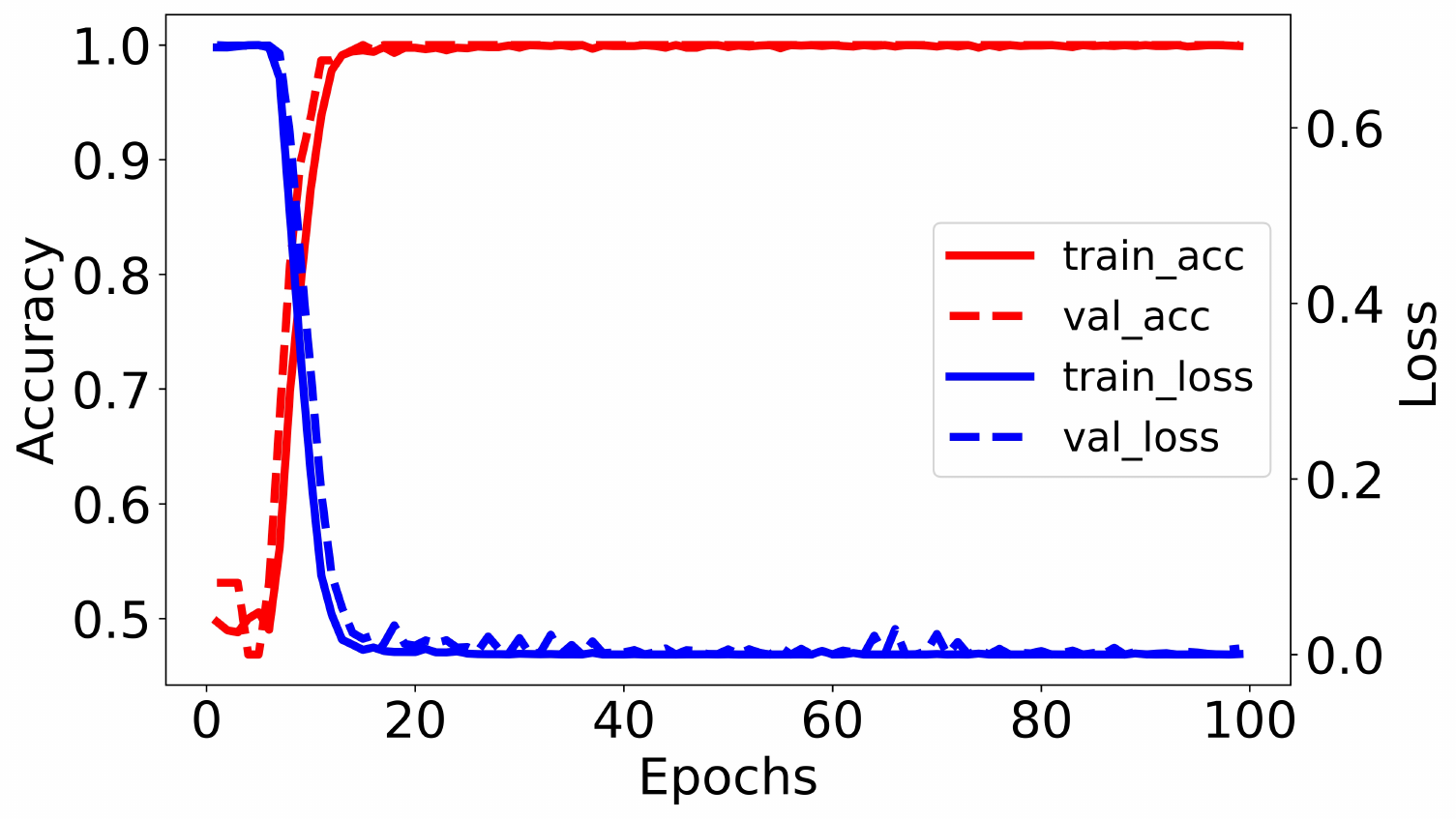}}
    \caption{Accuracy and loss curves during training.  The first, second, and third rows are the results of one-, two-, and three-dimensional convolutional neural networks in two different data sets, respectively. For each subgraph, the abscissa is the number of iterations, and the left and right ordinates are the accuracy and loss values, respectively.}
	\label{fig:acc loss curves}
\end{figure}

In this section, we report a series of experimental results aimed at comparing the classification performance of the one-, two-, and three-dimensional convolutional neural network in the diagnosis of PD. Tables \ref{tab:compare result on Dra} and Table \ref{tab:compare result on Pa} show the dynamic features used in the data enhancement process and the corresponding model diagnostic results obtained. To obtain robust experimental conclusions, both the {\em DraWritePD} data set \cite{Nomm_at_al_Luria_analysis_ICMLA2018} and the {\em PaHaW} data set \cite{drotar, drotar2014analysis} were considered. The related state-of-the-art results obtained on the {\em PaHaW} data set are shown in Table \ref{tab:discuss}. First, we evaluated the impact of embedding different sets of dynamic features in the same-dimensional space on the diagnostic performance of the model. Specifically, in the baseline experiment of each dimension space, only the coordinate feature or its derived velocity feature is used to encode the position information, and then, on this basis, other dynamic features gradually encode the colour information and line width information in the enhanced data. The experimental results in Table \ref{tab:compare result on Dra} and Table \ref{tab:compare result on Pa} show that, in the same dimensional space, overall diagnostic performance predictably presents the same upward trend, and the continuous improvement of performance confirms that encoding more dynamic features in the enhancement of the data helps to distinguish PD patients from HC subjects. In particular, the encoding of colour information and line width information makes the enhanced data more discriminating. There is, however, one exception. In the one-dimensional space, compared with the baseline experiment, the addition of velocity features failed to provide reasonable predictions, possibly because the velocity features were derived from coordinate features, resulting in redundant discriminating information. Furthermore, we compared the diagnostic performance of convolutional neural networks in different dimensions. It is worth noting that for a fair comparison, all convolutional neural network models used the same network architecture. First, the experimental results in Table \ref{tab:compare result on Dra} and Table \ref{tab:compare result on Pa} confirm that, with the same dynamic feature encoding, on average, convolutional neural networks achieve increasingly better diagnostic performance with increasing dimensionality, with the most competitive diagnostic results obtained by $3$D convolutional neural networks. For example, in terms of diagnostic accuracy, the $1$D, $2$D and $3$D convolutional networks achieved sequentially increasing diagnostic performance of $77.33\%$, $80.00\%$ and $86.67\%$, respectively, in the {\em PaHaW} data set. In addition, interestingly, the diagnostic performance of $3$D convolutional neural networks is almost comparable to that of low-dimensional convolutional neural networks even if only location information is given.  For example, in the DaWritePD data set, the $3$D convolutional neural network can achieve $77.78\%$ diagnostic accuracy only in the baseline experiment; on the contrary, the optimal diagnostic accuracy of $1$D and $2$D convolutional neural network $1$ D and $2$ D is just $62.5\%$ and $77.78\%$ respectively. Finally, the accuracy curves and loss curves in the specific training process are shown in Figure \ref{fig:acc loss curves}.

\begin{table}[htbp]
	\centering
	\footnotesize
    \caption{Performance comparison in the {\em DraWritePD} dataset.} 
    \label{tab:compare result on Dra}
    \begin{threeparttable}
    \setlength{\tabcolsep}{1.8mm}{
		\begin{tabular}{l l l l l l l l c c c c c}
			\toprule
			 & \multicolumn{6}{l}{Dynamic Features} & &\multicolumn{5}{l}{Metrics (in$\%$)} \\
            \cline{2-7} \cline{9-13} 
			Dimension & $x$ & $y$ & $a$ & $l$ & $p$ & $v$ & & Precision & Sensitivity & Specificity & Accuracy & $F_{1}$ score \\
            \hline
            $1$D & \ding{52} & \ding{52} &  &  &  &  &&  $50.00$ & $50.00$ & $50.00$ & $50.00$ & $50.00$ \\
               & \ding{52} & \ding{52} &  &  &  & \ding{52} &&  $50.00$ & $75.00$ & $25.00$ & $50.00$ & $60.00$ \\
               & \ding{52} & \ding{52} & \ding{52} & \ding{52} & \ding{52} &  &&  
            $50.00$ & $75.00$ &  $25.00$ & $50.00$ & $60.00$ \\
            & \ding{52} & \ding{52} & \ding{52} & \ding{52} & \ding{52} & \ding{52} && $60.00$  & $75.00$ &  $50.00$ & $62.50$ & $66.67$ \\
            \hline
            $2$D & \ding{52} & \ding{52} &  &  &  &  && $50.00$ & $50.00$ & $60.00$ & $55.56$ & $50.00$  \\
               & \ding{52} & \ding{52} &  &  &  & \ding{52} && $66.67$ & $50.00$ & $80.00$ & $66.67$ &  $57.14$ \\
               & \ding{52} & \ding{52} & \ding{52} & \ding{52} & \ding{52} &  &&  $60.00$ & $75.00$ &  $60.00$  &  $66.67$  &   $66.67$  \\
               & \ding{52} & \ding{52} & \ding{52} & \ding{52} & \ding{52} & \ding{52} && $75.00$ &  
                $75.00$ & $80.00$ & $77.78$ & $75.00$ \\
               \hline
            $3$D & \ding{52} & \ding{52} &  &  &  & \ding{52} && $75.00$ & $75.00$ & $80.00$ & $77.78$ &  $75.00$ \\
               & \ding{52} & \ding{52} & \ding{52} & \ding{52} & \ding{52} & \ding{52} && $77.50$ & $87.50$ & $80.00$ & $83.34$ & $81.95$  \\
			\toprule
		\end{tabular} }
        \begin{tablenotes}
        \footnotesize
        \item \textit{Note}: The abbreviations $x$, $y$ denote the x- and y- coordinate features, and $a$, $l$, $p$ and $v$ are the azimuth, altitude, pressure and velocity features, respectively. 
	\end{tablenotes}
 \end{threeparttable}
\end{table}

\begin{table}[htbp]
	\centering
	\footnotesize
    \caption{Performance comparison in the {\em PaHaW} dataset.} 
    \label{tab:compare result on Pa}
    \begin{threeparttable}
    \setlength{\tabcolsep}{1.8mm}{
		\begin{tabular}{l l l l l l l l c c c c c}
			\toprule
			 & \multicolumn{6}{l}{Dynamic Features} & &\multicolumn{5}{l}{Metrics (in$\%$)} \\
            \cline{2-7} \cline{9-13} 
			Dimension & $x$ & $y$ & $a$ & $l$ & $p$ & $v$ & & Precision & Sensitivity & Specificity & Accuracy & $F_{1}$ score \\
            \hline
            $1$D & \ding{52} & \ding{52} &  &  &  &  &&   $50.00$ &  
                 $57.14$ &  $50.00$ & $53.33$ &  $53.33$ \\
               & \ding{52} & \ding{52} &  &  &  & \ding{52} &&   $50.00$ &  $57.14$ & $50.00$ &   $53.33$ &  $53.33$  \\
               & \ding{52} & \ding{52} & \ding{52} & \ding{52} & \ding{52} &  &&  
             $57.14$ &  $57.14$ & $62.50$ & $60.00$ &  $57.14$ \\
            & \ding{52} & \ding{52} & \ding{52} & \ding{52} & \ding{52} & \ding{52} &&  $66.67$  &  $85.71$ &  $62.50$ & $73.33$ &  $75.00$\\
            \hline
            $2$D & \ding{52} & \ding{52} &  &  &  &  && $50.00$ & $85.71$ & $25.00$ & $53.33$ & $63.16$ \\
               & \ding{52} & \ding{52} &  &  &  & \ding{52} &&  $71.43$ & $71.43$ & $75.00$ & $73.33$ & $71.43$\\
               & \ding{52} & \ding{52} & \ding{52} & \ding{52} & \ding{52} &  && $83.33$ & $71.43$ & $87.50$ & $80.00$ &  $76.92$   \\
               & \ding{52} & \ding{52} & \ding{52} & \ding{52} & \ding{52} & \ding{52} &&  $75.00$ & $85.71$ & $75.00$ & $80.00$ & $80.00$  \\
               \hline
            $3$D & \ding{52} & \ding{52} &  &  &  & \ding{52} && $60.00$ & $85.71$ & $50.00$ & $66.67$ &  $70.59$ \\
               & \ding{52} & \ding{52} & \ding{52} & \ding{52} & \ding{52} & \ding{52} &&  $85.71$ &  $85.71$ & $87.50$ & $86.67$ &  $85.71$  \\
			\toprule
		\end{tabular} }
        \begin{tablenotes}
        \footnotesize
        \item \textit{Note}: The abbreviations $x$, $y$ denote the x- and y- coordinate features, and $a$, $l$, $p$ and $v$ are the azimuth, altitude, pressure and velocity features, respectively. 
	\end{tablenotes}
 \end{threeparttable}
\end{table}

\begin{table*}[htbp]
	\centering
	\footnotesize
    \caption{Comparisons with state-of-the-art works.} 
    \label{tab:discuss}
		\begin{tabular}{l l l l c}
			\toprule
			Author(s) & Dataset &  Features & Models & Accuracy (in $\%$) \\
            \hline 
              Drotár et al. \cite{drotar2016evaluation} & PaHaW & hand-crafted & SVM & $62.80$  \\
              Diaz et al. \cite{diaz2021sequence} & PaHaW & 1D CNN-extracted & 1D CNN + BiGRU &  $93.75$ \\
              Diaz et al. \cite{diaz2019dynamically} & PaHaW & 2D CNN-extracted & 2D CNN + SVM & $75.00$ \\
              Present work & PaHaW & 1D CNN-extracted & 1D CNN & $73.33$ \\
               &  & 2D CNN-extracted & 2D CNN & $80.00$ \\
               &  & 3D CNN-extracted & 3D CNN & $86.67$ \\
              
			\toprule
		\end{tabular} 
\end{table*}

\section{Discussion}\label{sec:discussion}

Comparing Tables \ref{tab:compare result on Dra} \ref{tab:compare result on Pa} that summarise model goodness metrics (computed on the basis of testing data) for the different feature sets and model dimensionalities one can observe that encoding more features into the image usually causes model goodness to increase and increasing dimensionality of the convolutional kernel also leads to better models. In the case of the {\em DraWritePD} data set, exceptions occur with specificity and precision. In the case of the {\em PaHaW} data set, exceptions also occur in the sensitivity of the models. Such anomalies in the behaviour of goodness metrics may be caused by the presence and absence of the feature describing the velocity that is related to the amount of tremor in the motions. Currently, there are many new features (spiral specific features) or improvements to existing features \cite{valla2022tremor} that can further improve the classification performance, similar to the diagnostic accuracy in $90\%$ obtained using multiple raw and derived features in \cite{diaz2021sequence}. Although such an investigation is beyond the scope of this paper, it may constitute the direction of future studies. Another point to discuss is that while the differences in performance metrics are similar between one- and two-dimensional CNNs, the difference between two- and three-dimensional CNNs is greater for the {\em DraWritePD} dataset. This may be triggered by the fact that in the case of {\em PaHaW} no reference drawing is provided, but in the case of {\em DraWritePD} a reference drawing is presented, making it easier to complete the test. The direction of the drawing and the age groups of the subjects tested may also contribute to this difference.


\section{Conclusions}\label{sec:conclusiona}
The application of one-, two-, and three-dimensional deep convolutional neural networks for the analysis of spiral drawing tests to support the diagnosis of Parkinson's disease has been investigated. Through comparative experiments on the two datasets, our hypothesis is confirmed that data representation and classification models in high-dimensional space are more beneficial to distinguish PD patients from HC subjects. Additionally, although we adapted the raw dynamic feature set for feature encoding to obtain high diagnostic accuracy, we believe that there is still room for improvement. Future research will be directed towards determining the specific feature sets to be encoded. 

\section*{Ethics statement}
The data acquisition process was carried out with strict guidance from the privacy law. The study was approved by the Research Ethics Committee of the University of Tartu (No. $1275T-9$). 

\section*{Author contributions}
S. N{\~o}mm problem conceptualisation, workflow curating, manuscript editing. X. Wang Workflow design, computational experiments, and manuscript writing. J. Huang and M. Chatzakou curating of the computational experiments, manuscript curating. 
E. Valla conceptualisation of the 3D data representation. Kadri Medijainen has performed data acquisition. Pille Taba, A. Toomela medical protocols curating. M. Ruzhansky manuscript, workflow, computations curating, supervision.

\section{Acknowledgements}
The authors thank P. Drotar and his team for sharing the PaHaW dataset. 
This work was partially supported by the FWO Odysseus 1 grant G.0H94.18N: Analysis and Partial Differential Equations and the Methusalem Programme of the Ghent University Special Research Fund (BOF) (Grant number 01M01021). Michael Ruzhansky is also supported by the EPSRC grant EP/R003025/2. Marianna Chatzakou is a postdoctoral fellow of the Research Foundation – Flanders (FWO)  under the postdoctoral grant No 12B1223N. The work of Sven N{\~o}mm and Elli Valla in the project ``ICT programme'' was supported by the European Union through the European Social Fund. Pille Taba's work was supported by the PRG grant $957$ from the Estonian Research Council.
\bibliographystyle{elsarticle-num} 
\bibliography{references}
\end{document}